\documentclass{article}
\usepackage{amssymb}
\usepackage{amsmath}

\input{tcilatex}
\begin{document}

\begin{center}
{\LARGE Simultaneous Diagonalization}

{\LARGE and SVD of\ Commuting Matrices\medskip }

{\large Ronald P. Nordgren\footnote{%
email: nordgren@rice.edu}}

Brown School of Engineering, Rice University{\LARGE \medskip }
\end{center}

\textbf{Abstract. }We present a matrix version of a known method of
constructing common eigenvectors of two diagonalizable commuting matrices,
thus enabling their simultaneous diagonalization. The matrices may have
simple eigenvalues of multiplicity greater than one. The singular value
decomposition (SVD) of a class of commuting matrices also is treated. The
effect of row/column permutation is examined. Examples are given.

\section{Introduction}

It is well known that if two diagonalizable matrices have the same
eigenvectors, then they commute. The converse also is true and a
construction for the common eigenvectors (enabling simultaneous
diagonalization) is known. If one of the matrices has distinct eigenvalues
(multiplicity one), it is easy to show that its eigenvectors pertain to both
commuting matrices. The case of matrices with simple eigenvalues of
multiplicity greater than one requires a more complicated construction of
their common eigenvectors. Here, we present a matrix version of a
construction procedure given by Horn and Johnson \cite[Theorem 1.3.12]{HORN}
and in a video by Sadun \cite{SADU}. The eigenvector construction procedure
also is applied to the singular value decomposition of a class of commuting
matrices that includes the case where at least one of the matrices is real
and symmetric. In addition, we consider row/column permutation of the
commuting matrices. Three examples illustrate the eigenvector construction
procedure.

\section{Eigenvector Construction}

Let $\mathbf{A}$ and $\mathbf{B}$ be diagonalizable square matrices that
commute, i.e.%
\begin{equation}
\mathbf{AB}=\mathbf{BA}  \label{ABBA}
\end{equation}%
and their Jordan canonical forms read%
\begin{equation}
\mathbf{A}=\mathbf{S}_{A}\mathbf{D}_{A}\mathbf{S}_{A}^{-1},\quad \mathbf{B}=%
\mathbf{S}_{B}\mathbf{D}_{B}\mathbf{S}_{B}^{-1}.  \label{JordAB}
\end{equation}%
Here, the columns of $\mathbf{S}_{A}$ are the eigenvectors $\mathbf{s}_{i}$
of $\mathbf{A}$ corresponding to the eigenvalues $\lambda _{i}$ in the
diagonal matrix $\mathbf{D}_{A}$, i.e.%
\begin{equation}
\mathbf{As}_{i}=\lambda _{i}\mathbf{s}_{i}.  \label{Asi}
\end{equation}%
If $\mathbf{A}$ is normal $\left( \mathbf{AA}^{\ast }=\mathbf{A}^{\ast }%
\mathbf{A}\right) $, then $\mathbf{S}_{A}$ can be made unitary $\left( 
\mathbf{S}_{A}^{-1}=\mathbf{S}_{A}^{\ast }\right) $. From (\ref{Asi}) and (%
\ref{ABBA}) we have%
\begin{equation}
\mathbf{BAs}_{i}=\mathbf{A}\left( \mathbf{Bs}_{i}\right) =\lambda _{i}\left( 
\mathbf{Bs}_{i}\right) ,  \label{ABsi}
\end{equation}%
whence $\mathbf{Bs}_{i}$ also is an eigenvector of $\mathbf{A}$ for the
eigenvalue $\lambda _{i}.$

When $\lambda _{i}$ are distinct, it follows from (\ref{Asi}) and (\ref{ABsi}%
) that each eigenvector $\mathbf{Bs}_{i}$ must be a scalar multiple of $%
\mathbf{s}_{i}$ and, noting (\ref{JordAB}), we have%
\begin{equation}
\mathbf{BS}_{A}=\mathbf{S}_{A}\mathbf{D}_{B}\quad \text{or}\quad \mathbf{D}%
_{B}=\mathbf{S}_{A}^{-1}\mathbf{BS}_{A},\quad \mathbf{B}=\mathbf{S}_{A}%
\mathbf{D}_{B}\mathbf{S}_{A}^{-1}  \label{BSD}
\end{equation}%
which shows that $\mathbf{S}_{A}$ is an eigenvector matrix for $\mathbf{B}$
as well as $\mathbf{A}$ and thus $\mathbf{S}_{A}$ diagonalizes both of them.

In the case that $\mathbf{A}$ has $k$ simple eigenvalues of multiplicity
greater than or equal to one, \textbf{$D$}$_{A}$ can be written in block
matrix form as%
\begin{equation}
\mathbf{D}_{A}=\limfunc{diag}\left[ \mathbf{D}_{A1},\mathbf{D}_{A2},\ldots ,%
\mathbf{D}_{Ak}\right] ,\quad \mathbf{D}_{Ai}=\limfunc{diag}\left[ \lambda
_{i},\lambda _{i},\ldots ,\lambda _{i}\right] =\lambda _{i}\mathbf{I}_{i},
\label{DA}
\end{equation}%
where $\mathbf{I}_{i}$ is the identity matrix of order equal to the number
of $\lambda _{i}$'s. Then $\mathbf{S}_{A}$ has the block form%
\begin{equation}
\mathbf{S}_{A}=\left[ \mathbf{S}_{A1},\mathbf{S}_{A2},\ldots ,\mathbf{S}_{Ak}%
\right] ,  \label{SA}
\end{equation}%
where the columns of $\mathbf{S}_{Ai}$ are the eigenvectors corresponding to 
$\lambda _{i}.$ In view of (\ref{ABsi}), the columns of $\mathbf{BS}_{Ai}$
also are eigenvectors corresponding to $\lambda _{i}$ and therefore they
must be linear combinations of the eigenvectors in $\mathbf{S}_{Ai},$ i.e.%
\begin{equation}
\mathbf{BS}_{Ai}=\mathbf{S}_{Ai}\mathbf{T}_{i},  \label{BSST}
\end{equation}%
where $\mathbf{T}_{i}$ is a square matrix of the same order as $\mathbf{I}%
_{i}.$ Thus, we may write $\mathbf{BS}_{A}$ as 
\begin{align}
\mathbf{BS}_{A}& =\left[ \mathbf{BS}_{A1},\mathbf{BS}_{A2},\ldots ,\mathbf{BS%
}_{Ak}\right] =\left[ \mathbf{S}_{A1}\mathbf{T}_{1},\mathbf{S}_{A2}\mathbf{T}%
_{2},\ldots ,\mathbf{S}_{Ak}\mathbf{T}_{k}\right]   \notag \\
& =\mathbf{S}_{A}\mathbf{T},\quad \mathbf{T}=\limfunc{diag}\left[ \mathbf{T}%
_{1},\mathbf{T}_{2},\ldots ,\mathbf{T}_{k}\right] =\mathbf{S}_{A}^{-1}%
\mathbf{BS}_{A}.  \label{BSZ}
\end{align}%
Since $\mathbf{B}$ is diagonalizable, each $\mathbf{T}_{i}$ is
diagonalizable and its Jordan form reads%
\begin{equation}
\mathbf{T}_{i}=\mathbf{S}_{Ti}\mathbf{D}_{Ti}\mathbf{S}_{Ti}^{-1},
\label{ZiJord}
\end{equation}%
where we take $\mathbf{S}_{Ti}=\mathbf{I}_{i}$ if $\lambda _{i}=0$ in order
to make$\mathbf{\ S}_{Ti}^{-1}=\mathbf{I}_{i}$ Then $\mathbf{T}$ has the
Jordan form 
\begin{equation}
\mathbf{T}=\mathbf{S}_{T}\mathbf{D}_{T}\mathbf{S}_{T}^{-1},  \label{ZJord}
\end{equation}%
where%
\begin{equation}
\mathbf{S}_{T}=\limfunc{diag}\left[ \mathbf{S}_{T1},\mathbf{S}_{T2},\ldots ,%
\mathbf{S}_{Tk}\right] ,\quad \mathbf{D}_{T}=\limfunc{diag}\left[ \mathbf{D}%
_{T1},\mathbf{D}_{T2},\ldots ,\mathbf{D}_{Tk}\right] .  \label{SZDZ}
\end{equation}%
From (\ref{BSZ}) and (\ref{ZJord}) we have 
\begin{equation}
\mathbf{B}=\mathbf{S}_{A}\mathbf{TS}_{A}^{-1}=\left( \mathbf{S}_{A}\mathbf{S}%
_{T}\right) \mathbf{D}_{T}\left( \mathbf{S}_{A}\mathbf{S}_{T}\right) ^{-1}
\label{BJord2}
\end{equation}%
which is a Jordan form for $\mathbf{B.}$ It follows that $\mathbf{S}_{A}%
\mathbf{S}_{T}$ is an eigenvector matrix for $\mathbf{B,}$ and $\mathbf{D}%
_{T}$ must contain the same eigenvalues as $\mathbf{D}_{B},$ but they may be
in a different order as shown by the examples below. 

From (\ref{JordAB}) we form%
\begin{equation}
\mathbf{A}=\left( \mathbf{S}_{A}\mathbf{S}_{T}\right) \left( \mathbf{S}%
_{T}^{-1}\mathbf{D}_{A}\mathbf{S}_{T}\right) \left( \mathbf{S}_{A}\mathbf{S}%
_{T}\right) ^{-1}.  \label{ASZ}
\end{equation}%
By (\ref{SZDZ}) we find that%
\begin{align}
\mathbf{S}_{T}^{-1}\mathbf{D}_{A}\mathbf{S}_{T}& =\limfunc{diag}\left[ 
\mathbf{S}_{T1}^{-1},\mathbf{S}_{T2}^{-1},\ldots ,\mathbf{S}_{Tk}^{-1}\right]
\limfunc{diag}\left[ \mathbf{D}_{A1},\mathbf{D}_{A2},\ldots ,\mathbf{D}_{Ak}%
\right]   \notag \\
& \qquad \qquad \times \limfunc{diag}\left[ \mathbf{S}_{T1},\mathbf{S}%
_{T2},\ldots ,\mathbf{S}_{Tk}\right]   \notag \\
& =\limfunc{diag}\left[ \mathbf{S}_{T1}^{-1}\lambda _{1}\mathbf{I}_{1}%
\mathbf{S}_{T1},\mathbf{S}_{T2}^{-1}\lambda _{2}\mathbf{I}_{2}\mathbf{S}%
_{T2},\ldots ,\mathbf{S}_{Tk}^{-1}\lambda _{k}\mathbf{I}_{k}\mathbf{S}_{Tk}%
\right]   \label{SZDA} \\
& =\limfunc{diag}\left[ \lambda _{1}\mathbf{I}_{1},\lambda _{2}\mathbf{I}%
_{2},\ldots ,\lambda _{k}\mathbf{I}_{k}\right]   \notag \\
& =\limfunc{diag}\left[ \mathbf{D}_{A1},\mathbf{D}_{A2},\ldots ,\mathbf{D}%
_{Ak}\right] =\mathbf{D}_{A}  \notag
\end{align}%
and (\ref{ASZ}) becomes%
\begin{equation}
\mathbf{A}=\left( \mathbf{S}_{A}\mathbf{S}_{T}\right) \mathbf{D}_{A}\left( 
\mathbf{S}_{A}\mathbf{S}_{T}\right) ^{-1}  \label{AST}
\end{equation}%
which shows that $\mathbf{S}_{A}\mathbf{S}_{T}$ is an eigenvector matrix for 
$\mathbf{A}$ as well as $\mathbf{B}.$ Thus, $\mathbf{A}$ and $\mathbf{B}$
can be simultaneously diagonalized by $\mathbf{S}_{A}\mathbf{S}_{T}$, i.e.%
\begin{equation}
\mathbf{D}_{A}=\left( \mathbf{S}_{A}\mathbf{S}_{T}\right) ^{-1}\mathbf{A}%
\left( \mathbf{S}_{A}\mathbf{S}_{T}\right) ,\quad \mathbf{D}_{B}=\left( 
\mathbf{S}_{A}\mathbf{S}_{T}\right) ^{-1}\mathbf{B}\left( \mathbf{S}_{A}%
\mathbf{S}_{T}\right) .  \label{Diag}
\end{equation}

We note that the role of $\mathbf{A}$ and $\mathbf{B}$ can be interchanged
in the above construction process. However, this results in essentially the
same common eigenvalue matrix as $\mathbf{S}_{A}\mathbf{S}_{T}$. To see
this, we rewrite (\ref{Diag}) and (\ref{BJord2}) as%
\begin{equation}
\mathbf{A=}\left( \mathbf{S}_{A}\mathbf{S}_{TA}\right) \mathbf{D}_{A}\left( 
\mathbf{S}_{A}\mathbf{S}_{TA}\right) ^{-1},\quad \mathbf{B=}\left( \mathbf{S}%
_{A}\mathbf{S}_{T}\right) \mathbf{D}_{T}\left( \mathbf{S}_{A}\mathbf{S}%
_{T}\right) ^{-1},  \label{ABre}
\end{equation}%
where, as noted above, $\mathbf{D}_{T}$ contains the same eigenvalues as the
original $\mathbf{D}_{B}$ but in a different position on the diagonal. Thus%
\begin{equation}
\mathbf{D}_{T}=\mathbf{P}_{B}\mathbf{D}_{B}\mathbf{P}_{B}^{T},\quad \mathbf{%
B=}\left( \mathbf{S}_{A}\mathbf{S}_{TA}\mathbf{P}_{B}\right) \mathbf{D}%
_{B}\left( \mathbf{S}_{A}\mathbf{S}_{TA}\mathbf{P}_{B}\right) ^{-1},
\label{DBre}
\end{equation}%
where $\mathbf{P}_{B}$ is a permutation matrix. Similarly the construction
starting with $\mathbf{B}$ results in%
\begin{gather}
\mathbf{D}_{A}=\mathbf{P}_{A}\mathbf{D}_{A}\mathbf{P}_{A}^{T},\quad \mathbf{%
A=}\left( \mathbf{S}_{B}\mathbf{S}_{TB}\mathbf{P}_{A}\right) \mathbf{D}%
_{A}\left( \mathbf{S}_{B}\mathbf{S}_{TB}\mathbf{P}_{A}\right) ^{-1},  \notag
\\
\mathbf{B=}\left( \mathbf{S}_{B}\mathbf{S}_{TB}\right) \mathbf{D}_{B}\left( 
\mathbf{S}_{B}\mathbf{S}_{TB}\right) ^{-1}.  \label{BAre}
\end{gather}%
On comparing these results, we see that%
\begin{equation}
\mathbf{S}_{B}\mathbf{S}_{TB}=\mathbf{S}_{A}\mathbf{S}_{TA}\mathbf{P}%
_{B},\quad \mathbf{S}_{A}\mathbf{S}_{TA}=\mathbf{S}_{B}\mathbf{S}_{TB}%
\mathbf{P}_{A},\quad \mathbf{P}_{A}=\mathbf{P}_{B}^{T},  \label{SbSa}
\end{equation}%
i.e. $\mathbf{S}_{B}\mathbf{S}_{TB}$ is a reordering of the eigenvectors
(columns) of $\mathbf{S}_{A}\mathbf{S}_{TA}$ according to the reordering of
the eigenvalues in $\mathbf{D}_{T}$ via $\mathbf{P}_{B}$ as an example below
illustrates.

Furthermore, if $\mathbf{A}$ has distinct eigenvalues, by (\ref{BSZ}) and (%
\ref{BSD}) (\ref{ZJord}) we have%
\begin{gather}
\mathbf{T}=\mathbf{S}_{A}^{-1}\mathbf{BS}_{A}=\mathbf{S}_{A}^{-1}\mathbf{S}%
_{A}\mathbf{D}_{B}\mathbf{S}_{A}^{-1}\mathbf{S}_{A}=\mathbf{ID}_{T}\mathbf{I}%
=\mathbf{S}_{T}\mathbf{D}_{T}\mathbf{S}_{T}^{-1},  \notag \\
\therefore \mathbf{S}_{T}=\mathbf{I},\quad \mathbf{S}_{A}\mathbf{S}_{T}=%
\mathbf{S}_{A},  \label{TSA}
\end{gather}%
so nothing is gained from the construction of $\mathbf{S}_{A}\mathbf{S}_{T}$
in this case.

\section{Singular Value Decomposition}

The construction of matrices in the singular value decompositions (SVD) of
two commuting matrices also is of interest. The SVD of $\mathbf{A}$ and $%
\mathbf{B}$ read%
\begin{equation}
\mathbf{A=U}_{A}\mathbf{\Sigma }_{A}\mathbf{V}_{A}^{\ast },\quad \mathbf{B=U}%
_{B}\mathbf{\Sigma }_{B}\mathbf{V}_{B}^{\ast },  \label{ABsvd}
\end{equation}%
where $\mathbf{U}_{A},\ \mathbf{V}_{A},\ \mathbf{U}_{B},$ and $\mathbf{V}_{B}
$ are unitary matrices, whereas $\mathbf{\Sigma }_{A}$ and $\mathbf{\Sigma }%
_{B}$ are diagonal matrices with non-negative real numbers (called singular
values) on the diagonal. It follows from (\ref{ABsvd}) that%
\begin{equation}
\mathbf{AA}^{\ast }\mathbf{=U}_{A}\mathbf{\Sigma }_{A}^{2}\mathbf{U}%
_{A}^{\ast },\quad \mathbf{BB}^{\ast }\mathbf{=U}_{B}\mathbf{\Sigma }_{B}^{2}%
\mathbf{U}_{B}^{\ast }  \label{AAt}
\end{equation}%
which are Jordan forms of $\mathbf{AA}^{\ast }$ and $\mathbf{BB}^{\ast }.$
These Jordan forms enable the determination of $\mathbf{U}_{A},\ \mathbf{U}%
_{B},\ \mathbf{\Sigma }_{A},$ and $\mathbf{\Sigma }_{B}$ after which $%
\mathbf{V}_{A}$ and $\mathbf{V}_{B}$ can be determined from (\ref{ABsvd}).
If $\mathbf{\Sigma }_{A}$ and $\mathbf{\Sigma }_{B}$ are nonsingular, then (%
\ref{ABsvd}) leads to%
\begin{equation}
\mathbf{V}_{A}=\mathbf{A}^{\ast }\mathbf{U}_{A}\mathbf{\Sigma }%
_{A}^{-1},\quad \mathbf{V}_{B}=\mathbf{B}^{\ast }\mathbf{U}_{B}\mathbf{%
\Sigma }_{B}^{-1}.  \label{VA}
\end{equation}%
and if they are singular, an alternate approach is given by Meyer \cite{MEYE}%
.

If $\mathbf{A\ }$and $\mathbf{B}$ commute and if 
\begin{equation}
\mathbf{A}^{\ast }\mathbf{B=BA}^{\ast }  \label{A*B}
\end{equation}%
then%
\begin{equation}
\mathbf{AA}^{\ast }\mathbf{BB}^{\ast }=\mathbf{ABA}^{\ast }\mathbf{B}^{\ast
}=\mathbf{BAB}^{\ast }\mathbf{A}^{\ast }=\mathbf{BB}^{\ast }\mathbf{AA}%
^{\ast },  \label{AABB}
\end{equation}%
i.e. $\mathbf{AA}^{\ast }$ and $\mathbf{BB}^{\ast }$ commute. Thus, they
have a common left-singular vector matrix $\mathbf{U}_{A}=\mathbf{U}_{B}$
which can be found by the foregoing eigenvector construction procedure. If $%
\mathbf{A\ }$and $\mathbf{B}$ also are normal with common unitary matrix $%
\mathbf{S}$, then $\mathbf{U}_{A}=\mathbf{U}_{B}=\mathbf{S.}$ Note that (\ref%
{A*B}) is satisfied if $\mathbf{A}$ is real and symmetric. The construction
of SVD matrices is illustrated in the examples below.

\section{Permutation}

We can form two new matrices by row/column permutations of the restricted
form%
\begin{equation}
\mathbf{\hat{A}=PAP}^{T},\quad \mathbf{\hat{B}=PBP}^{T},  \label{ABhat}
\end{equation}%
where $\mathbf{P}$ is a permutation matrix. Since $\mathbf{P}$ is orthogonal 
$\left( \mathbf{PP}^{T}=\mathbf{I}\right) ,$ it follows that $\mathbf{\hat{A}%
}$ and $\mathbf{\hat{B}}$ commute when $\mathbf{A}$ and $\mathbf{B}$ commute
as seen from 
\begin{equation}
\mathbf{\hat{A}\hat{B}=PAP}^{T}\mathbf{PBP}^{T}=\mathbf{PABP}^{T}=\mathbf{%
PBAP}^{T}=\mathbf{PB\mathbf{P}^{T}\mathbf{P}AP}^{T}=\mathbf{\hat{B}\hat{A}.}
\end{equation}%
Furthermore, the Jordan forms of $\mathbf{\hat{A}}$ and $\mathbf{\hat{B},}$
with (\ref{JordAB}), read 
\begin{gather}
\mathbf{\hat{A}=\mathbf{\hat{S}}}_{A}\mathbf{\mathbf{D}}_{A}\mathbf{\mathbf{%
\hat{S}}}_{A}^{-1},\quad \mathbf{\hat{B}=\mathbf{\hat{S}}}_{B}\mathbf{%
\mathbf{D}}_{B}\mathbf{\mathbf{\hat{S}}}_{B}^{-1},  \notag \\
\mathbf{\mathbf{\hat{S}}}_{A}=\mathbf{P\mathbf{S}}_{A},\quad \mathbf{\mathbf{%
\hat{S}}}_{B}=\mathbf{P\mathbf{S}}_{B}  \label{Jord ABh}
\end{gather}%
which show that the eigenvalues are unchanged by the permutation (\ref{ABhat}%
) and the eigenvectors are permuted. Similar formulas apply to the SVD's of $%
\mathbf{\hat{A}}$ and $\mathbf{\hat{B}}$ and their singular values also are
unchanged by the permutation (\ref{ABhat}). Indeed, it is known that the
singular values of any matrix $\mathbf{A}$ are invariant under row/column
permutations of the more general form 
\begin{equation}
\mathbf{\tilde{A}=PAQ,}  \label{PAQ}
\end{equation}%
where $\mathbf{Q}$ is a second permutation matrix. Next, we present three
examples to illustrate our theoretical results.

\section{Examples}

\noindent \textbf{Example 1}. We start with the normal matrix $\mathbf{A}$
and the magic square matrix $\mathbf{B}$ given by%
\begin{equation}
\mathbf{A}=\left[ 
\begin{array}{ccc}
1+i & 1 & 1 \\ 
1 & 1+i & 1 \\ 
1 & 1 & 1+i%
\end{array}%
\right] ,\quad \mathbf{B}=\left[ 
\begin{array}{ccc}
7 & 0 & 5 \\ 
2 & 4 & 6 \\ 
3 & 8 & 1%
\end{array}%
\right]   \label{AB33}
\end{equation}%
which commute and whose Jordan-form matrices are%
\begin{align}
\mathbf{S}_{A}& =\frac{1}{6}\left[ 
\begin{array}{ccc}
2\sqrt{3} & 3\sqrt{2} & \sqrt{6} \\ 
2\sqrt{3} & 0 & -2\sqrt{6} \\ 
2\sqrt{3} & -3\sqrt{2} & \sqrt{6}%
\end{array}%
\right] ,\quad \mathbf{S}_{B}=\left[ 
\begin{array}{ccc}
1 & 5 & 5 \\ 
1 & 2\left( 1+\sqrt{6}\right)  & 2\left( 1-\sqrt{6}\right)  \\ 
1 & -\left( 7+2\sqrt{6}\right)  & -\left( 7-2\sqrt{6}\right) 
\end{array}%
\right] ,  \notag \\
\mathbf{D}_{A}& =\limfunc{diag}\left[ 3+i,i,i\right] ,\quad \mathbf{D}_{B}=%
\limfunc{diag}\left[ 12,-2\sqrt{6},2\sqrt{6}\right] ,  \label{SD3}
\end{align}
where $\mathbf{S}_{A}$ is orthogonal. Since $\mathbf{B}$ has distinct
eigenvalues, $\mathbf{S}_{B}$ is a (nonorthogonal) eigenvector matrix for $%
\mathbf{A}$ as well as $\mathbf{B}$. However, $\mathbf{S}_{A}$ is not an
eigenvector matrix for $\mathbf{B}$ since $\mathbf{A}$ has multiple
eigenvalues.

In the SVD matrices for $\mathbf{A}$ and $\mathbf{B,}$ noting that (\ref{A*B}%
) is satisfied, we find that $\mathbf{U}_{A}=\mathbf{U}_{B}=\mathbf{S}_{A}$
and 
\begin{align}
\mathbf{\Sigma }_{A}& =\limfunc{diag}\left[ \sqrt{10},1,1\right] ,\quad 
\mathbf{V}_{A}=\frac{1}{30}\left[ 
\begin{array}{ccc}
\sqrt{30}\left( 3-i\right)  & -15i\sqrt{2} & -5i\sqrt{6} \\ 
\sqrt{30}\left( 3-i\right)  & 0 & 10i\sqrt{6} \\ 
\sqrt{30}\left( 3-i\right)  & 15i\sqrt{2} & -5i\sqrt{6}%
\end{array}%
\right] ,  \notag \\
\mathbf{\Sigma }_{B}& =\limfunc{diag}\left[ 12,4\sqrt{3},2\sqrt{3}\right]
,\quad \mathbf{V}_{B}=\frac{1}{6}\left[ 
\begin{array}{ccc}
2\sqrt{3} & \sqrt{6} & 3\sqrt{2} \\ 
2\sqrt{3} & -2\sqrt{6} & 0 \\ 
2\sqrt{3} & \sqrt{6} & -3\sqrt{2}%
\end{array}%
\right] .  \label{VAB3}
\end{align}

\noindent \textbf{Example 2}. The example given in the video by Sadun \cite%
{SADU} has the commuting matrices%
\begin{equation}
\mathbf{A}=\left[ 
\begin{array}{cccc}
0 & 4 & 0 & 0 \\ 
1 & 0 & 0 & 0 \\ 
0 & 0 & 0 & 4 \\ 
0 & 0 & 1 & 0%
\end{array}%
\right] ,\quad \mathbf{B}=\left[ 
\begin{array}{cccc}
0 & 0 & 1 & 0 \\ 
0 & 0 & 0 & 1 \\ 
1 & 0 & 0 & 0 \\ 
0 & 1 & 0 & 0%
\end{array}%
\right]  \label{AB4}
\end{equation}
whose Jordan-form matrices are%
\begin{align}
\mathbf{S}_{A}& =\left[ 
\begin{array}{cccc}
2 & 0 & 2 & 0 \\ 
1 & 0 & -1 & 0 \\ 
0 & 2 & 0 & 2 \\ 
0 & 1 & 0 & -1%
\end{array}%
\right] ,\quad \mathbf{S}_{B}=\frac{\sqrt{2}}{2}\left[ 
\begin{array}{cccc}
0 & 1 & 0 & 1 \\ 
1 & 0 & 1 & 0 \\ 
0 & 1 & 0 & -1 \\ 
1 & 0 & -1 & 0%
\end{array}%
\right] ,  \label{SD4} \\
\mathbf{D}_{A}& =\limfunc{diag}\left[ 2,2,-2,-2\right] ,\quad \mathbf{D}_{B}=%
\limfunc{diag}\left[ 1,1,-1,-1\right] ,  \notag
\end{align}%
where $\mathbf{S}_{B}$ is orthogonal.

On following the matrix construction procedure for the common eigenvector
matrix from $\mathbf{A},$ we find that%
\begin{gather}
\mathbf{T}_{A}=\mathbf{S}_{A}^{-1}\mathbf{BS}_{A}=\left[ 
\begin{array}{cccc}
0 & 1 & 0 & 0 \\ 
1 & 0 & 0 & 0 \\ 
0 & 0 & 0 & 1 \\ 
0 & 0 & 1 & 0%
\end{array}%
\right] ,\quad \mathbf{S}_{TA}=\left[ 
\begin{array}{cccc}
1 & 1 & 0 & 0 \\ 
-1 & 1 & 0 & 0 \\ 
0 & 0 & 1 & 1 \\ 
0 & 0 & -1 & 1%
\end{array}%
\right] ,  \notag \\
\mathbf{S}_{A}\mathbf{S}_{TA}=\left[ 
\begin{array}{cccc}
2 & 2 & 2 & 2 \\ 
1 & 1 & -1 & -1 \\ 
-2 & 2 & -2 & 2 \\ 
-1 & 1 & 1 & -1%
\end{array}%
\right] ,\quad 
\begin{array}{c}
\mathbf{D}_{A}=\limfunc{diag}\left[ 2,2,-2,-2\right] , \\ 
\mathbf{\tilde{D}}_{B}=\limfunc{diag}\left[ -1,1,-1,1\right] .%
\end{array}
\label{SAZ4}
\end{gather}%
It can be verified that $\mathbf{S}_{A}\mathbf{S}_{TA}$ is an eigenvalue
matrix for both $\mathbf{A}$ and $\mathbf{B}.$ Note that the eigenvalues in $%
\mathbf{\tilde{D}}_{B}$ are in a different order than those in $\mathbf{D}%
_{B}$ and they are related by (\ref{DBre}) with $\mathbf{D}_{T}\equiv 
\mathbf{\tilde{D}}_{B}$ and 
\begin{equation}
\mathbf{P}_{B}=\left[ 
\begin{array}{cccc}
0 & 0 & 0 & 1 \\ 
0 & 1 & 0 & 0 \\ 
0 & 0 & 1 & 0 \\ 
1 & 0 & 0 & 0%
\end{array}%
\right] .  \label{PB}
\end{equation}

On following the matrix construction procedure for the common eigenvector
matrix from $\mathbf{B}$ (instead of $\mathbf{A})$ we find that%
\begin{gather}
\mathbf{T}_{B}=\mathbf{S}_{B}^{-1}\mathbf{AS}_{B}=\left[ 
\begin{array}{cccc}
0 & 1 & 0 & 0 \\ 
4 & 0 & 0 & 0 \\ 
0 & 0 & 0 & 1 \\ 
0 & 0 & 4 & 0%
\end{array}%
\right] ,\quad \mathbf{S}_{TB}=\left[ 
\begin{array}{cccc}
-1 & 1 & 0 & 0 \\ 
2 & 2 & 0 & 0 \\ 
0 & 0 & -1 & 1 \\ 
0 & 0 & 2 & 2%
\end{array}%
\right] ,  \notag \\
\mathbf{S}_{B}\mathbf{S}_{TB}=\left[ 
\begin{array}{cccc}
2 & 2 & 2 & 2 \\ 
-1 & 1 & -1 & 1 \\ 
2 & 2 & -2 & -2 \\ 
-1 & 1 & 1 & -1%
\end{array}%
\right] ,\quad 
\begin{array}{c}
\mathbf{\tilde{D}}_{A}=\limfunc{diag}\left[ -2,2,-2,2\right] , \\ 
\mathbf{D}_{B}=\limfunc{diag}\left[ 1,1,-1,-1\right] ,%
\end{array}%
,  \label{SBZ4}
\end{gather}%
and (\ref{SbSa}) can be verified.

Noting that $\mathbf{AA}^{\ast }$ and $\mathbf{BB}^{\ast }$ commute since $%
\mathbf{B}$ is real and symmetric, the SVD matrices for $\mathbf{A}$ and $%
\mathbf{B}$ are%
\begin{align}
\mathbf{U}_{A}& \mathbf{=}\mathbf{\mathbf{U}}_{B}\mathbf{=I,\quad V}_{A}=%
\left[ 
\begin{array}{cccc}
0 & 1 & 0 & 0 \\ 
1 & 0 & 0 & 0 \\ 
0 & 0 & 0 & 1 \\ 
0 & 0 & 1 & 0%
\end{array}%
\right] \mathbf{,}  \label{UV4A} \\
\mathbf{\Sigma }_{A}& =\limfunc{diag}\left[ 4,1,4,1\right] ,\quad \mathbf{%
\Sigma }_{B}=\mathbf{I,\quad \mathbf{V}}_{B}\mathbf{=B.}  \notag
\end{align}%
An alternate set of SVD matrices for $\mathbf{A}$ and $\mathbf{B}$ is%
\begin{align}
\mathbf{U}_{A}& \mathbf{=\mathbf{U}}_{B}\mathbf{=S}_{B},\quad \mathbf{\Sigma 
}_{A}=\limfunc{diag}\left[ 1,4,1,4\right] ,\quad \mathbf{\Sigma }_{B}=%
\mathbf{I,}  \notag \\
\mathbf{V}_{A}& \mathbf{=}\frac{\sqrt{2}}{2}\left[ 
\begin{array}{cccc}
1 & 0 & 1 & 0 \\ 
0 & 1 & 0 & 1 \\ 
1 & 0 & -1 & 0 \\ 
0 & 1 & 0 & -1%
\end{array}%
\right] ,\quad \mathbf{\mathbf{\mathbf{V}}}_{B}\mathbf{=}\frac{\sqrt{2}}{2}%
\left[ 
\begin{array}{cccc}
0 & 1 & 0 & -1 \\ 
1 & 0 & -1 & 0 \\ 
0 & 1 & 0 & 1 \\ 
1 & 0 & 1 & 0%
\end{array}%
\right] ,  \label{UV4B}
\end{align}%
where $\mathbf{\mathbf{\mathbf{V}}}_{A}$ is determined from (\ref{VA}) and $%
\mathbf{\mathbf{\mathbf{V}}}_{B}$ is formed by changing the sign of the two
eigenvectors (last two columns) in $\mathbf{S}_{B}$ associated with its
negative eigenvalues $\left( -1,-1\right) $ in (\ref{SD4}).

\noindent \textbf{Example 3. }Following a method given by Nordgren \cite%
{NORD}, we construct the symmetric commuting matrices%
\begin{equation}
\mathbf{A}=\left[ 
\begin{array}{cccccc}
1 & 0 & 2 & 3 & 0 & 4 \\ 
0 & 3 & 0 & 0 & 7 & 0 \\ 
2 & 0 & 1 & 4 & 0 & 3 \\ 
3 & 0 & 4 & 1 & 0 & 2 \\ 
0 & 7 & 0 & 0 & 3 & 0 \\ 
4 & 0 & 3 & 2 & 0 & 1%
\end{array}%
\right] ,\quad \mathbf{B}=\left[ 
\begin{array}{cccccc}
i & i & i & 1 & 1 & 1 \\ 
i & i & i & 1 & 1 & 1 \\ 
i & i & i & 1 & 1 & 1 \\ 
1 & 1 & 1 & i & i & i \\ 
1 & 1 & 1 & i & i & i \\ 
1 & 1 & 1 & i & i & i%
\end{array}%
\right] ,  \label{AB6}
\end{equation}%
with orthogonal Jordan-form matrices%
\begin{align}
\mathbf{S}_{A}& =\frac{1}{4}\left[ 
\begin{array}{cccccc}
-2 & \sqrt{2} & \sqrt{2} & -2 & \sqrt{2} & \sqrt{2} \\ 
0 & 2 & -2 & 0 & 2 & -2 \\ 
2 & \sqrt{2} & \sqrt{2} & 2 & \sqrt{2} & \sqrt{2} \\ 
2 & -\sqrt{2} & -\sqrt{2} & -2 & \sqrt{2} & \sqrt{2} \\ 
0 & -2 & 2 & 0 & 2 & -2 \\ 
-2 & -\sqrt{2} & -\sqrt{2} & 2 & \sqrt{2} & \sqrt{2}%
\end{array}%
\right] , \\
\mathbf{D}_{A}& =\limfunc{diag}\left[ 0,-4,-4,-2,10,10\right] ,  \notag
\end{align}%
\begin{align}
\mathbf{S}_{B}& =\frac{\sqrt{6}}{12}\left[ 
\begin{array}{cccccc}
2 & 2 & 2 & 2 & 2 & 2 \\ 
2 & -\sqrt{3}-1 & \sqrt{3}-1 & 2 & -\sqrt{3}-1 & \sqrt{3}-1 \\ 
2 & \sqrt{3}-1 & -\sqrt{3}-1 & 2 & \sqrt{3}-1 & -\sqrt{3}-1 \\ 
-2 & -2 & -2 & 2 & 2 & 2 \\ 
-2 & \sqrt{3}+1 & -\sqrt{3}+1 & 2 & -\sqrt{3}-1 & \sqrt{3}-1 \\ 
-2 & -\sqrt{3}+1 & \sqrt{3}+1 & 2 & \sqrt{3}-1 & -\sqrt{3}-1%
\end{array}%
\right] , \\
\mathbf{D}_{B}& =\limfunc{diag}\left[ -3+3i,0,0,3+3i,0,0\right] .  \notag
\end{align}%
The matrix construction procedure leads to the following common eigenvector
matrix and corresponding eigenvalues:%
\begin{align}
\mathbf{S}_{A}\mathbf{S}_{TA}& =\frac{\sqrt{3}}{6}\left[ 
\begin{array}{cccccc}
-\sqrt{3} & \sqrt{2} & 1 & -\sqrt{3} & \sqrt{2} & 1 \\ 
0 & \sqrt{2} & -2 & 0 & \sqrt{2} & -2 \\ 
\sqrt{3} & \sqrt{2} & 1 & \sqrt{3} & \sqrt{2} & 1 \\ 
\sqrt{3} & -\sqrt{2} & -1 & -\sqrt{3} & \sqrt{2} & 1 \\ 
0 & -\sqrt{2} & 2 & 0 & \sqrt{2} & -2 \\ 
-\sqrt{3} & -\sqrt{2} & -1 & \sqrt{3} & \sqrt{2} & 1%
\end{array}%
\right] ,  \label{SAST6} \\
\mathbf{D}_{A}& =\limfunc{diag}\left[ -2,-2,-2,0,8,8,\right] ,\quad \mathbf{D%
}_{B}=\limfunc{diag}\left[ 0,-3+3i,0,0,3+3i,0\right] .  \notag
\end{align}

Since $\mathbf{A}$ and $\mathbf{B}$ are symmetric and $\mathbf{S}_{A}\mathbf{%
S}_{TA}$ is orthogonal, suitable SVD matrices are%
\begin{align}
\mathbf{U}_{A}& =\mathbf{U}_{B}=\mathbf{S}_{A}\mathbf{S}_{TA},  \notag \\
\mathbf{V}_{A}& =\frac{\sqrt{3}}{6}\left[ 
\begin{array}{cccccc}
\sqrt{3} & -\sqrt{2} & -1 & \sqrt{3} & \sqrt{2} & 1 \\ 
0 & -\sqrt{2} & 2 & 0 & \sqrt{2} & -2 \\ 
-\sqrt{3} & -\sqrt{2} & -1 & -\sqrt{3} & \sqrt{2} & 1 \\ 
-\sqrt{3} & \sqrt{2} & 1 & \sqrt{3} & \sqrt{2} & 1 \\ 
0 & \sqrt{2} & -2 & 0 & \sqrt{2} & -2 \\ 
\sqrt{3} & \sqrt{2} & 1 & -\sqrt{3} & \sqrt{2} & 1%
\end{array}%
\right] ,  \label{UV6} \\
\mathbf{V}_{B}& =\frac{\sqrt{3}}{6}\left[ 
\begin{array}{cccccc}
-\sqrt{3} & -1-i & 1 & -\sqrt{3} & 1-i & 1 \\ 
0 & -1-i & -2 & 0 & 1-i & -2 \\ 
\sqrt{3} & -1-i & 1 & \sqrt{3} & 1-i & 1 \\ 
\sqrt{3} & 1+i & -1 & -\sqrt{3} & 1-i & 1 \\ 
0 & 1+i & 2 & 0 & 1-i & -2 \\ 
-\sqrt{3} & 1+i & -1 & \sqrt{3} & 1-i & 1%
\end{array}%
\right] , \\
\mathbf{\Sigma }_{A}& =\limfunc{diag}\left[ 0,4,4,2,10,10\right] ,\quad 
\mathbf{\Sigma }_{B}=\limfunc{diag}\left[ 0,3\sqrt{2},0,0,3\sqrt{2},0\right]
,  \notag
\end{align}%
where $\mathbf{V}_{A}$ is formed by changing the sign of the first three
columns of $\mathbf{S}_{A}\mathbf{S}_{TA}$ associated with the negative
eigenvalues in $\mathbf{D}_{A}.$

To illustrate row/column permutation, let%
\begin{equation}
\mathbf{P}=\left[ 
\begin{array}{cccccc}
0 & 1 & 0 & 0 & 0 & 0 \\ 
0 & 0 & 0 & 0 & 0 & 1 \\ 
0 & 0 & 1 & 0 & 0 & 0 \\ 
1 & 0 & 0 & 0 & 0 & 0 \\ 
0 & 0 & 0 & 1 & 0 & 0 \\ 
0 & 0 & 0 & 0 & 1 & 0%
\end{array}%
\right] .
\end{equation}%
By (\ref{ABhat}), we have%
\begin{equation}
\mathbf{\hat{A}}=\left[ 
\begin{array}{cccccc}
3 & 0 & 0 & 0 & 0 & 7 \\ 
0 & 1 & 3 & 4 & 2 & 0 \\ 
0 & 3 & 1 & 2 & 4 & 0 \\ 
0 & 4 & 2 & 1 & 3 & 0 \\ 
0 & 2 & 4 & 3 & 1 & 0 \\ 
7 & 0 & 0 & 0 & 0 & 3%
\end{array}%
\right] ,\quad \mathbf{\hat{B}}=\left[ 
\begin{array}{cccccc}
i & 1 & i & i & 1 & 1 \\ 
1 & i & 1 & 1 & i & i \\ 
i & 1 & i & i & 1 & 1 \\ 
i & 1 & i & i & 1 & 1 \\ 
1 & i & 1 & 1 & i & i \\ 
1 & i & 1 & 1 & i & i%
\end{array}%
\right] .
\end{equation}%
It can be verified that $\mathbf{\hat{A}}$ and $\mathbf{\hat{B}}$ commute
and that their eigenvalues and singular values are the same as those of $%
\mathbf{A}$ and $\mathbf{B}$.

\section{Conclusion}

The three examples illustrate the efficacy of the presented matrix
construction procedure for the common eigenvectors of two commuting
matrices, thereby enabling their simultaneous diagonalization. The
construction also is useful in finding the SVD of matrices when at least one
of them is real and symmetric. A restricted row/column permutation of two
commuting matrices produces two computing matrices with unchanged
eigenvalues and singular values.

\end{document}